\documentstyle[12pt,amssymb]{article}

\begin{document}

\begin{center}
{\large {\bf A product formula of Mumford forms, and the rationality}} 
\vspace{1ex}

{\large {\bf of Ruelle zeta values for Schottky groups}} 
\end{center}

\begin{center}
{\sc Takashi Ichikawa} 
\end{center}
\vspace{0ex}

\begin{center}
{\bf Abstract} 
\end{center}

In any positive genus case, 
we show an explicit formula of the Mumford forms expressed by 
infinite products like Selberg type zeta values for Schottky groups. 
This result is considered as an extension of the formula in terms of 
Ramanujan's Delta function in the genus $1$ case, 
and is especially applied to studying the rationality 
of Ruelle zeta values for Schottky groups. 
\vspace{2ex}

\begin{center}
{\bf 1. Introduction} 
\end{center}

Let ${\mathcal M}_{g}$ be the moduli stack over ${\mathbb Z}$ of 
proper smooth (algebraic) curves of genus $g > 0$, 
and $\pi : {\mathcal C}_{g} \rightarrow {\mathcal M}_{g}$ be 
the universal curve. 
Then for each positive integer $k$, 
we have the $k$th tautological bundle, 
and the $k$th tautological line bundle 
$$
\Lambda_{k} = \pi_{*} 
\left( \Omega_{{\mathcal C}_{g}/{\mathcal M}_{g}}^{\otimes k} \right), 
\ \mbox{and} \ 
\lambda_{k} = \det \left( \Lambda_{k} \right)  
$$
on ${\mathcal M}_{g}$ respectively. 
Put $d_{k} = 6k^{2} - 6k + 1$. 
Then by a result of Mumford [M2], 
$\lambda_{k}$ is isomorphic to $\lambda_{1}^{\otimes d_{k}}$, 
and hence there exists a nowhere vanishing global section 
$\mu_{g;k}$ of 
$\lambda_{k}^{\otimes (-1)} \otimes \lambda_{1}^{\otimes d_{k}}$ 
(with certain boundary condition) 
which is called the {\it Mumford form} 
and is uniquely determined up to a sign. 
It is shown in [D2] and [GS] using Arakelov theory that $\mu_{g;k}$ 
gives an isometry (up to a constant) between $\lambda_{k}$ and 
$\lambda_{1}^{\otimes d_{k}}$ with Quillen metric. 

In particular, $\mu_{g;2}$ was well studied since $\lambda_{2}$ is 
isomorphic to the volume form on ${\mathcal M}_{g}$ by Kodaira-Spencer theory, 
and $\mu_{g;2}$ is connected with bosonic string measure. 
As is well known, 
under the trivialization of 
$\lambda_{2}^{\otimes (-1)} \otimes \lambda_{1}^{\otimes 13}$ arising from 
the natural homomorphism ${\rm Sym}^{2}(\Lambda_{1}) \rightarrow \Lambda_{2}$, 
$\mu_{1;2}$ is expressed by Ramanujan's Delta function $\Delta(\tau)$ as 
$$
\mu_{1;2} = \Delta(\tau) \stackrel{\rm def}{=} 
q \prod_{m=1}^{\infty} \left( 1 - q^{m} \right)^{24}; \ 
q = e^{2 \pi \sqrt{-1} \tau} \ \left( {\rm Im}(\tau) > 0 \right), 
$$
and $\mu_{2;2} = \left( \theta_{2}/2^{6} \right)^{2}$, 
where $\theta_{g}$ denotes the Siegel modular form defined as 
the product of even theta constants of degree $g$. 
Furthermore, 
it is shown in [I2, 4] that $\mu_{3;2} = \sqrt{- \theta_{3} / 2^{28}}$, 
and that $\mu_{4;2}$ is expressed as certain derivatives of 
the normalized Schottky's Siegel modular form. 

In order to study Mumford forms for any $g$ and $k$, 
we consider Schottky uniformization theory which gives rise to 
local coordinates on the moduli space 
${\mathcal M}_{g}({\mathbb C})$ of Riemann surfaces of genus $g$ 
in terms of Schottky groups 
$\Gamma \subset PGL_{2}({\mathbb C})$ of rank $g$. 
Let $X_{\Gamma}$ denote the Riemann surface which is Schottky uniformized by 
$\Gamma$. 
Then McIntyre-Takhtajan [MT] constructed the period map which is 
a non-degenerate pairing between the space 
$H^{0} \left( X_{\Gamma}, \Omega^{k}_{X_{\Gamma}} \right)$ 
of holomorphic $k$-forms on $X_{\Gamma}$, 
and the Eichler cohomology group of $\Gamma$ with coefficients 
in polynomials of degree $\leq 2k - 2$. 
By fixing generators of $\Gamma$, 
we have the {\it normalized} basis of 
$H^{0} \left( X_{\Gamma}, \Omega^{k}_{X_{\Gamma}} \right)$ 
which is defined to be dual to a canonical ${\mathbb Z}$-basis of 
the Eichler cohomology group, 
and hence there is a local trivialization of $\lambda_{k}$ 
given by their exterior products. 
Under this trivialization, 
Zograf [Z1, 2] and McIntyre-Takhtajan [MT] gave the holomorphic factorization 
formula which represents an isomorphism between 
$\lambda_{k}$ with Quillen metric and the classical Liouville action 
in terms of the higher genus version of $\Delta(\tau)$. 
From their formula, 
we have the following expression of $\mu_{g;k}$ by an infinite product for 
$\Gamma = \langle \gamma_{1},..., \gamma_{g} \rangle$, 
namely, 
there exists a nonzero constant $c(g;k)$ depending only on $g$ and $k$ 
such that 
$$
c(g;k) \cdot \mu_{g;k} = 
\frac{\prod_{\{ \gamma \}} \left( 
\prod_{m=1}^{k-1} \left( 1 - q_{\gamma}^{m} \right)^{d_{k}} \cdot 
\prod_{m=k}^{\infty} \left( 1 - q_{\gamma}^{m} \right)^{d_{k}-1} \right)}
{\left( 1 - q_{\gamma_{1}} \right)^{2} \cdots 
\left( 1 - q_{\gamma_{1}}^{k-1} \right)^{2} 
\left( 1 - q_{\gamma_{2}}^{k-1} \right)}, 
$$
where $q_{\gamma}$ denotes the multiplier of $\gamma$ and 
$\{ \gamma\}$ runs over primitive conjugacy classes of $\Gamma$. 
Then a main result of this paper states: 
\vspace{2ex}

{\bf Theorem} (cf. Theorem 6.1). 
\begin{it} 
$c(g;k) = \pm 1$. 
\end{it}
\vspace{2ex}

Note that the right hand side of the above formula becomes 
$\prod_{m=1}^{\infty} (1 - q^{m})^{24}$ when $g = 1$, $k = 2$, 
and hence this theorem gives a generalization of the product formula 
of $\mu_{1;2}$ using $\Delta(\tau)$ 
(the absence of the factor $q$ comes from the difference of trivializations 
of $\lambda_{2}$). 
This theorem is also extended in Theorem 3.4 for pointed algebraic curves. 

The proof of this theorem is obtained by applying 
the arithmetic version of Schottky uniformization theory [I2] 
to the following simple fact in number theory: 
{\it if two integral power series are known to be proportional 
to each other 
and primitive (i.e., not congruent to $0$ modulo any prime), 
then these power series are equal up to a sign.} 
Here this uniformization theory gives rise to 
generalized Tate curves whose normalized regular $1$-forms 
$\left\{ \omega_{i} \right\}$ and multiplicative periods 
have universal expressions as integral power series. 
Since $\mu_{g;k}$ is a ``primitive'' object 
for the ${\mathbb Z}$-structure of ${\cal M}_{g}$, 
the key point of the proof is to show the assertion that 
the above trivialization of $\lambda_{k}$ comes from 
a basis of regular $k$-forms on a generalized Tate curve. 
In the previous version of the paper, 
we proved this assertion only for $k = 2$ by constructing a subset of 
$\left\{ \omega_{i} \omega_{j} \right\}$ which spans generically 
the normalized basis of 
$H^{0} \left( X_{\Gamma}, \Omega^{k}_{X_{\Gamma}} \right)$ 
over a field of any characteristic (see 6.2). 
In this revised version, 
we prove this assertion for any $k$ by studying the behavior of 
$\lambda_{k}$ and its trivialization associated with the period map 
under the degeneration and resolution of Riemann surfaces. 
Then our theorem is shown by considering such processes 
from a generalized Tate curve to a pointed Tate curve for which 
the normalized basis is easily constructed. 

As an application of the above theorem, 
we consider the rationality of the special values of the Ruelle zeta function 
$$
R_{\Gamma}(s) = \prod_{\{ \gamma \}} \left( 1 - |q_{\gamma}|^{s} \right)^{-1} 
$$
for $\Gamma$ based on the fact that if $\Gamma \subset PSL_{2}({\mathbb R})$, 
then the infinite products of Zograf and McIntyre-Takhtajan 
become Selberg type zeta values. 
A conjecture of Deligne [D1], 
which is extended by Beilinson [B] and Bloch-Kato [BK], 
states the following: 
each special value of the Hasse-Weil $L$-function of a motive ${\cal M}$ 
becomes the product of its transcendental part represented by 
the period or regulator of ${\cal M}$, 
and its rational part represented by the arithmetic invariant of ${\cal M}$. 
As its analogy for {\it geometric} zeta functions, 
we show the following: 
\vspace{2ex}

{\bf Theorem} (see Theorem 7.2 and Corollary 7.3 for the precise statement). 
\begin{it} 
Assume that $\Gamma \subset PSL_{2}({\mathbb R})$ and that the algebraic curve 
$X_{\Gamma}$ is defined over a subfield $K$ of ${\mathbb C}$. 
Let $c(\Gamma)$ denote the ``value'' at $\Gamma$ of Zograf's infinite product 
divided by the period of regular $1$-forms on $X_{\Gamma}$. 
Then the modified Ruelle zeta value of $\Gamma$ at any integer $k > 1$ 
belongs to 
$$
\frac{\mbox{period of regular $(k+1)$-forms on $X_{\Gamma}$}}
{\mbox{period of regular $k$-forms on $X_{\Gamma}$}} 
\cdot c(\Gamma)^{12k} \cdot K^{\times}. 
$$
\end{it} 
\vspace{-1ex}

We also express this zeta value by the discriminant of $X_{\Gamma}$ 
based on results of Saito [S], 
however we have no result on the rational or transcendental property of 
$c(\Gamma)$. 
Notice that using a result of [I3], 
one can construct Schottky groups $\Gamma \subset PSL_{2}({\mathbb R})$ 
such that $X_{\Gamma}$ are hyperelliptic and defined over ${\mathbb Q}$. 

By using the Selberg trace formula, 
the leading term at the central point $0$ of the Ruelle $L$ function 
of a hyperbolic manifold $M$ was studied by Fried [Fri] when $M$ is compact, 
and by Park [P], Sugiyama [Su1-5] and Gon-Park [GP] 
when $M$ has finite volume. 
Especially, 
it seems interesting that Sugiyama's results for hyperbolic $3$-folds 
are regarded as a geometric analog to Iwasawa theory 
and the Beilinson conjecture. 
In our case where $X_{\Gamma}$ is the boundary of the hyperbolic $3$-fold 
$M = {\mathbb H}^{3} / \Gamma$, 
the above theorem states that 
the rationality of the Ruelle zeta values of $M$ at positive integers 
is controlled by the field of definition of its boundary. 
\vspace{2ex}

\begin{center}
{\bf Acknowledgments} 
\end{center}

The author would like to thank deeply J. Park for his introduction 
of Zograf's formula which gave a motivation of this work. 
The author also would like to thank A. McIntyre, K. Sugiyama and S. Koyama 
for their valuable comments on this work. 
\vspace{2ex}

\begin{center}
{\bf 2. Formulas of Zograf and of McIntyre-Takhtajan} 
\end{center}

2.1. \ 
A {\it marked} Schottky group 
$\left( \Gamma; \gamma_{1},..., \gamma_{g} \right)$ 
consists of a Schottky group $\Gamma \subset PGL_{2}({\mathbb C})$ 
with its generators $\gamma_{1},..., \gamma_{g}$. 
Then taking the quotient by $\Gamma$ of its region of discontinuity 
$\Omega_{\Gamma} \subset {\mathbb P}_{\mathbb C}^{1}$, 
we have the Riemann surface $X_{\Gamma} = \Omega_{\Gamma} / \Gamma$ 
of genus $g$ which is called Schottky uniformized by $\Gamma$. 
For each marked Schottky group 
$\left( \Gamma; \gamma_{1},..., \gamma_{g} \right)$, 
there are open domains $D_{\pm i} \subset {\mathbb P}_{\mathbb C}^{1}$ 
$(1 \leq i \leq g)$ bounded by Jordan curves $C_{\pm i}$ such that 
$\gamma_{i} \left( {\mathbb P}_{\mathbb C}^{1} - D_{-i} \right)$ is 
the closure of $D_{i}$. 
Then 
$$
{\mathbb P}_{\mathbb C}^{1} - \left( \bigcup_{i=1}^{g} 
\left( D_{i} \cup D_{-i} \right) \right) 
$$
is a fundamental domain of $\Omega_{\Gamma}$ for $\Gamma$, 
and there is a basis $\omega_{1},..., \omega_{g}$ of holomorphic $1$-forms 
on $X_{\Gamma}$ which is {\it normalized} in the sense that 
${\displaystyle \frac{1}{2 \pi \sqrt{-1}} \oint_{C_{i}} \omega_{j}}$ 
is Kronecker's delta $\delta_{ij}$, 
where $C_{i}$ is counterclockwise oriented. 
Since $\Gamma - \{ 1 \}$ consists of hyperbolic, i.e., loxodromic elements, 
each $\gamma \in \Gamma - \{ 1 \}$ has uniquely 
{\it attractive} (resp. {\it repulsive}) fixed points 
$\alpha_{\gamma}$ (resp. $\beta_{\gamma}$) in ${\mathbb P}_{\mathbb C}^{1}$ 
and {\it multiplier} $q_{\gamma} \in {\mathbb C}$ satisfying 
$0 < |q_{\gamma}| < 1$, 
and 
$$
\frac{\gamma(z) - \alpha_{\gamma}}{\gamma(z) - \beta_{\gamma}} = 
q_{\gamma} \frac{z - \alpha_{\gamma}}{z - \beta_{\gamma}} 
\ \left( z \in {\mathbb P}_{\mathbb C}^{1} \right) 
$$
which means 
$$
\gamma = 
\left( \begin{array}{cc} \alpha_{\gamma} & \beta_{\gamma} \\ 1 & 1 
\end{array} \right) 
\left( \begin{array}{cc} 1     & 0      \\ 0 & q_{\gamma} \end{array} \right) 
\left( \begin{array}{cc} \alpha_{\gamma} & \beta_{\gamma} \\ 1 & 1 
\end{array} \right)^{-1} 
\ {\rm mod} \left( {\mathbb C}^{\times} \right). 
$$
For $\gamma = \gamma_{i}$, we put 
$\alpha_{\gamma} = \alpha_{i}$, $\beta_{\gamma} = \alpha_{-i}$ 
and $q_{\gamma} = q_{i}$, 
and we call $\left( \Gamma; \gamma_{1},..., \gamma_{g} \right)$ 
{\it normalized} 
if $\alpha_{1} = 0$, $\alpha_{-1} = \infty$ and $\alpha_{2} = 1$. 
Each marked Schottky group is conjugate to a unique normalized one, 
and if $q_{1},..., q_{g}$ are sufficiently small, 
then the exponent of convergence $\delta(\Gamma)$ of $\Gamma$ satisfies 
$\delta(\Gamma) < 1$. 

Let ${\mathfrak M}_{g}$ be the moduli space of 
compact Riemann surfaces $X$ of genus $g$ 
which is the complex orbifold associated with ${\mathcal M}_{g}$. 
Then $H^{0} \left( X, \Omega_{X}^{k} \right)$ gives 
the fiber of $\Lambda_{k}$ at the point corresponding to $X$. 
Let ${\mathfrak S}_{g}$ be the Schottky space of degree $g$ which is, 
by definition, 
a connected complex manifold classifying marked normalized Schottky groups 
$(\Gamma; \gamma_{1},..., \gamma_{g})$ of rank $g$. 
Then by $(\Gamma; \gamma_{1},..., \gamma_{g}) \mapsto X_{\Gamma}$, 
we have the natural covering map 
${\mathfrak S}_{g} \rightarrow {\mathfrak M}_{g}$, 
and the above $\omega_{1},..., \omega_{g}$ give a basis of $\Lambda_{1}$ 
over ${\mathfrak S}_{g}$. 

Let $S : {\mathfrak S}_{g} \rightarrow {\mathbb R}$ denote 
the classical Liouville action defined in [ZT]. 
\vspace{2ex}

2.2. \ 
We review results of Zograf [Z1, 2] and of McIntyre-Takhtajan [MT]. 
\vspace{2ex}

{\bf Zograf's formula} (cf. [Z1, 2] and [MT, Theorem 1]). 
\begin{it} 

{\rm (1)} 
There exists a holomorphic function $F(1)$ on ${\mathfrak S}_{g}$ which gives 
an isometry between $\lambda_{1}$ on ${\mathfrak M}_{g}$ with Quillen metric 
and the holomorphic line bundle on ${\mathfrak M}_{g}$ determined by 
the hermitian metric $\exp \left(S/12 \pi \right)$. 
More precisely, 
there exists a positive real number $c_{g}$ depending only on $g$ such that 
$$
\exp \left( \frac{S}{12 \pi} \right) = c_{g} \left| F(1) \right|^{2} 
\left\| \omega_{1} \wedge \cdots \wedge \omega_{g} \right\|_{\rm Q}^{2}, 
$$
where $\| * \|_{\rm Q}$ denotes the Quillen metric. 

{\rm (2)} 
On the open subset of ${\mathfrak S}_{g}$ corresponding to 
Schottky groups $\Gamma$ with $\delta(\Gamma) < 1$, 
$F(1)$ is given by the absolutely convergent infinite product  
$$
\prod_{\{ \gamma \}} \prod_{m=0}^{\infty} 
\left( 1 - q_{\gamma}^{1+m} \right), 
$$
where $\{ \gamma \}$ runs over primitive conjugacy classes of $\Gamma$. 
\end{it}
\vspace{2ex}

Assume that $g > 1$, and take an integer $k > 1$. 
Let $\left( \Gamma; \gamma_{1},..., \gamma_{g} \right)$ be 
a marked normalized Schottky group, 
and ${\mathbb C}[z]_{2k-2}$ be the ${\mathbb C}$-vector space of polynomials 
$f = f(z)$ of $z$ with degree $\leq 2k-2$ on which $\Gamma$ acts as 
$$
\gamma(f)(z) = f \left( \gamma(z) \right) \cdot \gamma'(z)^{1-k} \ 
\left( \gamma \in \Gamma, \ f \in {\mathbb C}[z]_{2k-2} \right). 
$$
Take $\xi_{1,k-1}$, $\xi_{2,1},..., \xi_{2,2k-2}$, 
$\xi_{i,0},..., \xi_{i,2k-2}$ $(3 \leq i \leq g)$ 
as elements of the Eichler cohomology group 
$H^{1} \left( \Gamma, {\mathbb C}[z]_{2k-2} \right)$ 
of $\Gamma$ which are uniquely determined by the condition: 
$$
\xi_{i,j}(\gamma_{l}) = \left\{ \begin{array}{ll} 
\delta_{2l} (z-1)^{j} & (i = 2), 
\\ 
\delta_{il} z^{j} & (i \neq 2) 
\end{array} \right. 
$$
for $1 \leq l \leq g$. 
Then it is shown in [MT, Section 4] that 
$$
\Psi_{g;k} (\varphi, \xi) \stackrel{\rm def}{=} 
\frac{1}{2 \pi \sqrt{-1}} \sum_{i=1}^{g} 
\oint_{C_{i}} \varphi \cdot \xi(\gamma_{i}) dz 
$$ 
for $\varphi \in H^{0} \left( X_{\Gamma}, \Omega^{k}_{X_{\Gamma}} \right)$, 
$\xi \in H^{1} \left( \Gamma, {\mathbb C}[z]_{2k-2} \right)$ 
is a non-degenerate pairing on 
$$
H^{0} \left( X_{\Gamma}, \Omega^{k}_{X_{\Gamma}} \right) \times 
H^{1} \left( \Gamma, {\mathbb C}[z]_{2k-2} \right). 
$$
Then there exists a basis 
$$ 
\left\{ \varphi_{1,k-1}, \ \varphi_{2,1},..., \varphi_{2,2k-2}, \ 
\varphi_{i,0},..., \varphi_{i,2k-2} \ (3 \leq i \leq g) \right\} 
$$ 
which we call the {\it normalized basis} of 
$H^{0} \left( X_{\Gamma}, \Omega^{k}_{X_{\Gamma}} \right)$ such that 
$\Psi_{g;k} \left( \varphi_{i,j}, \xi_{l,m} \right) = 
\delta_{il} \cdot \delta_{jm}$. 
\vspace{2ex} 

{\it Remark.} 
Since $-\pi \cdot \Psi_{g;k}$ is the pairing given in [MT, (4.1)], 
$$ 
\left\{ -\frac{\varphi_{1,k-1}}{\pi}, \ 
-\frac{\varphi_{2,1}}{\pi},..., -\frac{\varphi_{2,2k-2}}{\pi}, \ 
-\frac{\varphi_{j,0}}{\pi},..., -\frac{\varphi_{j,2k-2}}{\pi} \ 
(3 \leq j \leq g) \right\} 
$$ 
is the {\it natural basis for $n$-differentials} defined in [MT]. 
\vspace{2ex} 

In what follows, put 
$$
\left\{ \varphi_{1},..., \varphi_{(2k-1)(g-1)} \right\} = 
\left\{ \varphi_{1,k-1}, \ \varphi_{2,1},..., \varphi_{2,2k-2}, \ 
\varphi_{j,0},..., \varphi_{j,2k-2} \ (3 \leq j \leq g) \right\}. 
$$
\vspace{-1ex} 

{\bf McIntyre-Takhtajan's formula} (cf. [MT, Theorem 2]). 
\begin{it} 
Assume that $k > 1$. 

{\rm (1)} 
There exists a holomorphic function $F(k)$ on ${\mathfrak S}_{g}$ 
which gives an isometry between $\lambda_{k}$ on ${\mathfrak M}_{g}$ 
with Quillen metric and the holomorphic line bundle on ${\mathfrak M}_{g}$ 
determined by the hermitian metric $\exp \left( S/12 \pi \right)^{d_{k}}$. 
More precisely, 
there exists a positive real number $c_{g;k}$ depending only on $g$ and $k$ 
such that 
$$
\exp \left( \frac{S}{12 \pi} \right)^{d_{k}} = c_{g;k} \left| F(k) \right|^{2} 
\left\| \varphi_{1} \wedge \cdots \wedge \varphi_{(2k-1)(g-1)} 
\right\|_{\rm Q}^{2}, 
$$
where $\| * \|_{\rm Q}$ denotes the Quillen metric. 

{\rm (2)} 
On the whole Schottky space ${\mathfrak S}_{g}$ classifying 
marked normalized Schottky groups 
$\left( \Gamma; \gamma_{1},..., \gamma_{g} \right)$, 
$F(k)$ is given by the absolutely convergent infinite product 
$$
\left( 1 - q_{\gamma_{1}} \right)^{2} \cdots 
\left( 1 - q_{\gamma_{1}}^{k-1} \right)^{2} 
\left( 1 - q_{\gamma_{2}}^{k-1} \right) 
\prod_{\{ \gamma \}} \prod_{m=0}^{\infty} \left( 1 - q_{\gamma}^{k+m} \right), 
$$
where $\{ \gamma \}$ runs over primitive conjugacy classes of $\Gamma$. 
\end{it}
\vspace{2ex}

{\it Remark.} 
It is noticed in [MT, Remark 6] that by the extra factors of $F(k)$ in (2), 
the right hand side of the formula in (1) is invariant under permutations of 
the generators $\gamma_{1},..., \gamma_{g}$. 
\vspace{2ex}

2.3. \ 
From the above formulas, we have: 
\vspace{2ex}

{\bf Proposition 2.1.} \begin{it} 
There exists a nonzero constant $c(g;k)$ depending only on $g$ and $k$ 
such that 
$$
c(g;k) \cdot \mu_{g;k} \left( \varphi_{1} \wedge \cdots \wedge 
\varphi_{(2k-1)(g-1)} \right) 
= \frac{F(1)^{d_{k}}}{F(k)} 
\left( \omega_{1} \wedge \cdots \wedge \omega_{g} \right)^{\otimes d_{k}}, 
$$
and that 
$$
\left| c(g; k) \right|^{2} = \frac{c_{g;k}}{(c_{g})^{d_{k}}} 
\exp \left[ (g-1)(k^{2}-k) \left( 24 \zeta'_{\mathbb Q}(-1)-1 \right) \right], 
$$
where $\zeta'_{\mathbb Q}(s)$ denotes the derivative of 
the Riemann zeta function $\zeta_{\mathbb Q}(s)$. 
\end{it}
\vspace{2ex}

{\it Proof.} 
By the arithmetic Riemann-Roch theorem for algebraic curves (cf. [D2, GS]), 
the Mumford isomorphism $\mu_{g;k}$ gives rise to an isometry 
$$
(\lambda_{k})_{\rm Q} \cong 
(\lambda_{1})_{\rm Q}^{\otimes d_{k}} \cdot 
\exp \left[ a(g)(1 - d_{k})/12 \right] 
$$
between the metrized tautological line bundles with Quillen metric 
on ${\cal M}_{g}$, 
where $a(g) = (1-g)\left( 24 \zeta'_{\mathbb Q}(-1) -1 \right)$ 
is the Deligne constant. 
Therefore, by the formulas of Zograf and of McIntyre-Takhtajan, 
there exists a holomorphic function $c(g;k)$ on ${\mathfrak S}_{g}$ 
satisfying the above formula such that $\left| c(g;k) \right|$ is 
a constant function with value 
$$
\sqrt{c_{g;k}/(c_{g})^{d_{k}}} \exp \left[ a(g)(1 - d_{k})/12 \right]. 
$$
Since ${\mathfrak S}_{g}$ is connected, 
$c(g;k)$ is also a constant function. 
\ $\square$ 
\vspace{2ex}

\begin{center}
{\bf 3. Mumford isomorphism for pointed curves} 
\end{center}

3.1. \ 
Fix integers $g, n \geq 0$ such that $2g - 2 + n > 0$, 
and let $\overline{\cal M}_{g,n}$ be the moduli stack over ${\mathbb Z}$ 
of stable $n$-pointed curves of genus $g$. 
Then by definition, 
there exists the universal curve 
$\pi : \overline{\cal C}_{g,n} \rightarrow \overline{\cal M}_{g,n}$, 
and the universal sections 
$\sigma_{j} : \overline{\cal M}_{g,n} \rightarrow \overline{\cal C}_{g,n}$ 
$(1 \leq j \leq n)$. 
Denote by ${\cal M}_{g,n}$ the open substack of $\overline{\cal M}_{g,n}$ 
classifying proper smooth $n$-pointed curves of genus $g$, 
and put $\partial {\cal M}_{g,n} = \overline{\cal M}_{g,n} - {\cal M}_{g,n}$. 
Then we have the following line bundles on $\overline{\cal M}_{g,n}$: 
\begin{eqnarray*}
\lambda_{g,n;k} 
& \stackrel{\rm def}{=} & 
\det R \pi_{*} \left( 
\omega^{k}_{\overline{\cal C}_{g,n}/\overline{\cal M}_{g,n}} 
\left( k \sum_{j} \sigma_{j} \right) \right), 
\\ 
\psi_{g,n} 
& \stackrel{\rm def}{=} & 
\bigotimes_{j=1}^{n} \psi_{g,n}^{(j)}; \ 
\psi_{g,n}^{(j)} \stackrel{\rm def}{=} \sigma^{*}_{j} 
\left( \omega_{\overline{\cal C}_{g,n}/\overline{\cal M}_{g,n}} \right), 
\\ 
\delta_{g,n} 
& \stackrel{\rm def}{=} & 
{\cal O}_{\overline{\cal M}_{g,n}} \left( \partial {\cal M}_{g,n} \right). 
\end{eqnarray*} 
Using the $k$th residue map ${\rm Res}^{k}_{\sigma_{j}}$ given by 
${\rm Res}^{k}_{\sigma_{j}} \left( \eta \left( dz_{j}/z_{j} \right)^{k} 
\right) = \eta$, 
where $z_{j}$ denotes a local coordinate on 
$\overline{\cal C}_{g,n}$ around $\sigma_{j}$, 
we have 
$$
\lambda_{g,n;k} \cong 
\det \pi_{*} \left( 
\omega^{k}_{\overline{\cal C}_{g,n}/\overline{\cal M}_{g,n}} 
\left( (k - 1) \sum_{j} \sigma_{j} \right) \right). 
$$
In particular, 
$\lambda_{g,n;1} \cong \det \pi_{*} 
\left( \omega_{\overline{\cal C}_{g,n}/\overline{\cal M}_{g,n}} \right)$. 
When $n = 0$, 
we delete $n$ from these notations, 
and put $\psi_{g} = {\cal O}_{\overline{\cal M}_{g,n}}$. 
For a positive integer $k$, 
put $d_{k} = 6 k^{2} - 6 k + 1$. 
\vspace{2ex}

{\bf Theorem 3.1.} 
\begin{it} 
There exists a unique (up to a sign) isomorphism 
$$
\mu_{g,n;k} : \lambda_{g,n;k} \stackrel{\sim}{\rightarrow} 
\lambda_{g,n;1}^{\otimes d_{k}} \otimes 
\left( \psi_{g,n} \otimes \delta_{g,n}^{\otimes (-1)} 
\right)^{\otimes (k^{2}-k)/2}. 
$$
on $\overline{\cal M}_{g,n}$ which is called the Mumford isomorphism. 
\end{it} 
\vspace{2ex}

{\it Proof.} 
The uniqueness follows from the properness of $\overline{\cal M}_{g,n}$ 
over ${\mathbb Z}$. 
The existence is shown by Mumford [M2] when $n = 0$, 
namely we have an isomorphism 
$$
\lambda_{g+n;k} \cong \lambda_{g+n;1}^{\otimes d_{k}} \otimes 
\delta_{g+n}^{\otimes (k - k^{2})/2}. 
$$
Let 
${\rm cl}_{1} : 
\overline{\cal M}_{g,n} \times \overline{\cal M}_{1,1}^{\times n} 
\rightarrow \overline{\cal M}_{g+n}$ 
be the clutching morphism 
given by identifying the $j$th section $\sigma_{j}$ 
over $\overline{\cal M}_{g,n}$ and the unique section of the $j$th component 
of $\overline{\cal M}_{1,1}^{\times n}$ $(1 \leq j \leq n)$. 
Then by [Fre, Corollary 3.4], 
\begin{eqnarray*} 
{\rm cl}_{1}^{*} \left( \lambda_{g+n;k} \right) 
& \cong & 
\lambda_{g,n;k} \boxtimes \lambda_{1,1;k}^{\boxtimes n}, 
\\ 
{\rm cl}_{1}^{*} \left( \delta_{g+n} \right) 
& \cong & 
\left( \delta_{g,n} \otimes \psi_{g,n}^{\otimes (-1)} \right) \boxtimes 
\left( \delta_{1,1} \otimes \psi_{1,1}^{\otimes (-1)} \right)^{\boxtimes n}, 
\end{eqnarray*} 
and there exists an isomorphism between 
$\lambda_{g,n;k} \boxtimes \lambda_{1,1;k}^{\boxtimes n}$ and 
$$
\left( \lambda_{g,n;1}^{\otimes d_{k}} \otimes \left( 
\delta_{g,n} \otimes \psi_{g,n}^{\otimes (-1)} \right)^{\otimes (k - k^{2})/2} 
\right) \boxtimes 
\left( \lambda_{1,1;1}^{\otimes d_{k}} \otimes \left( 
\delta_{1,1} \otimes \psi_{1,1}^{\otimes (-1)} \right)^{\otimes (k - k^{2})/2} 
\right)^{\boxtimes n}. 
$$
Denote by 
$\Delta : {\rm Spec}({\mathbb Z}) \rightarrow 
\overline{\cal M}_{1,1}^{\times n}$ 
the section corresponding to the $n$ copies of the $1$-pointed stable curve 
obtained from ${\mathbb P}_{\mathbb Z}^{1}$ by identifying $0$ and $\infty$ 
with marked point $1$. 
Then taking the pullback of the above isomorphism by 
${\rm id} \times \Delta : \overline{\cal M}_{g,n} \rightarrow 
\overline{\cal M}_{g,n} \times \overline{\cal M}_{1,1}^{\times n}$, 
we have the required isomorphism. 
\ $\square$ 
\vspace{2ex} 

{\it Remark.} 
Let $\kappa_{g,n}$ be the line bundle on $\overline{\cal M}_{g,n}$ defined as 
the following Deligne's pairing: 
$$
\kappa_{g,n} = 
\left\langle 
\omega_{\overline{\cal C}_{g,n}/\overline{\cal M}_{g,n}} 
\left( \sum_{j} \sigma_{j} \right), 
\omega_{\overline{\cal C}_{g,n}/\overline{\cal M}_{g,n}} 
\left( \sum_{j} \sigma_{j} \right) 
\right\rangle. 
$$
Then by [W, 2.2] and [Fre, Theorem 3.5], 
$$
\lambda_{g,n;k}^{\otimes 12} \otimes \psi_{g,n} \otimes 
\delta_{g,n}^{\otimes (-1)} \cong \kappa_{g,n}^{\otimes d_{k}}, 
$$
and hence by comparing the isomorphism with that obtained by $k = 1$, 
we have 
$$
\lambda_{g,n;k}^{\otimes 12} \cong 
\lambda_{g,n;1}^{\otimes 12 d_{k}} \otimes 
\left( \psi_{g,n} \otimes \delta_{g,n}^{\otimes (-1)} \right)^{6(k^{2} - k)}. 
$$
This seems not imply directly the theorem since we do not know 
whether the Picard group of $\overline{\cal M}_{g,n}$ is torsion-free or not. 
\vspace{2ex}

3.2. \ 
Let 
${\rm cl}_{2} : \overline{\cal M}_{g,n+2} \rightarrow 
\overline{\cal M}_{g+1,n}$ 
be the clutching morphism given by identifying the sections 
$\sigma_{n+1}$ and $\sigma_{n+2}$. 
\vspace{2ex}

{\bf Proposition 3.2.} 
\begin{it} 
We have the following isomorphisms: 
\end{it}
\begin{eqnarray*}
{\rm cl}_{2}^{*} \left( \lambda_{g+1,n;k} \right) & \cong & 
\lambda_{g,n+2;k}, 
\\ 
{\rm cl}_{2}^{*} \left( \psi_{g+1,n} \right) & \cong & 
\psi_{g,n+2} \otimes \left( \psi_{g,n+2}^{(n+1)} \right)^{\otimes(-1)} 
\otimes \left( \psi_{g,n+2}^{(n+2)} \right)^{\otimes(-1)}, 
\\ 
{\rm cl}_{2}^{*} \left( \delta_{g+1,n} \right) & \cong & 
\delta_{g,n+2} \otimes \left( \psi_{g,n+2}^{(n+1)} \right)^{\otimes(-1)} 
\otimes \left( \psi_{g,n+2}^{(n+2)} \right)^{\otimes(-1)}. 
\end{eqnarray*}
Furthermore, 
${\rm cl}_{2}^{*} \left( \mu_{g+1,n;k} \right) = \pm \mu_{g,n+2;k}$. 
\vspace{2ex}

{\it Proof.} 
Let $(\pi : C \rightarrow S; \sigma_{1},..., \sigma_{n+2})$ 
be an $(n+2)$-pointed stable curve of genus $g$, 
and $(\pi' : C' \rightarrow S; \sigma_{1},..., \sigma_{n})$ 
be the $n$-pointed stable curve of $g + 1$ obtained from $C$ 
by identifying $\sigma_{n+1}$ and $\sigma_{n+2}$. 
Then as is shown in [K, Section 1], 
$\pi'_{*} \left( \omega^{k}_{C'} \left( k \sum_{j=1}^{n} \sigma_{j} 
\right) \right)$ 
is isomorphic to ${\rm Ker}(\varphi)$, 
where  
$\varphi : \pi_{*} \left( \omega^{k}_{C} \left( k \sum_{j=1}^{n+2} \sigma_{j} 
\right) \right) \rightarrow {\cal O}_{S}$ 
is given by 
$$
\varphi(\eta) = {\rm Res}^{k}_{\sigma_{n+1}} (\eta) - 
(-1)^{k} {\rm Res}^{k}_{\sigma_{n+2}} (\eta). 
$$
Therefore, we have the first isomorphism. 
The second and third ones are shown in [K, Theorem 4.3], 
and hence the last identity follows from 
the uniqueness of the Mumford isomorphism. 
\ $\square$ 
\vspace{2ex} 

3.3. \ 
Let $(\Gamma; \gamma_{1},..., \gamma_{g})$ 
be a marked Schottky group of rank $g > 1$, 
and $X_{\Gamma} = \Omega_{\Gamma}/\Gamma$ be the Riemann surface 
uniformized by $\Gamma$ with $n$-marked points 
given by $s_{1},..., s_{n} \in \Omega_{\Gamma}$. 
Denote by ${\mathbb C}[z]_{d}$ the ${\mathbb C}$-vector space of 
polynomials over ${\mathbb C}$ of $z$ with degree $\leq d$. 
For $k > 1$, 
we have a ${\mathbb C}$-bilinear form $\Psi_{g,n;k}$ on 
$$
H^{0} \left( X_{\Gamma}, 
\Omega^{k}_{X_{\Gamma}} \left( k \sum_{j} s_{j} \right) \right) 
\times 
\left( H^{1} \left(\Gamma, {\mathbb C}[z]_{2k-2} \right) \oplus 
\left( {\mathbb C}[z]_{k-1} \right)^{\oplus n} \right) 
$$  
which is defined as 
$$
\Psi_{g,n;k} \left( \varphi (dz)^{k}, \left( \xi, (f_{j})_{j} \right) \right) 
= \frac{1}{2 \pi \sqrt{-1}} \sum_{i=1}^{g} 
\oint_{\partial D_{i}} \varphi \cdot \xi(\gamma_{i}) dz 
+ \sum_{j=1}^{n} {\rm Res}_{s_{j}} (\varphi \cdot f_{j} dz). 
$$     

{\bf Proposition 3.3.} 
\begin{it} 
If $k > 1$, 
then $\Psi_{g,n;k}$ is non-degenerate. 
\end{it} 
\vspace{2ex}

{\it Proof.} 
This follows easily from that $\Psi_{g;k}$ in 2.2 is non-degenerate. 
\ $\square$ 
\vspace{2ex} 

Let 
$$
\left\{ \xi_{1,n-1}, \xi_{2,1},..., \xi_{2,2k-2}, 
\xi_{i,0},..., \xi_{i,2k-2} \ (3 \leq l \leq g) \right\} 
$$
be the basis of $H^{1} \left( \Gamma, {\mathbb C}[z]_{2k-2} \right)$ 
given in 2.2. 
Then this basis together with 
$$
\left\{ z^{d_{j}} \ (1 \leq j \leq n, 0 \leq d_{j} \leq k-1) \right\} 
$$
give a basis of 
$H^{1} \left( \Gamma, {\mathbb C}[z]_{2k-2} \right) \oplus 
\left( {\mathbb C}[z]_{k-1} \right)^{\oplus n}$. 
Put 
$$
h_{k} \stackrel{\rm def}{=} 
\dim H^{0} \left( X_{\Gamma}, 
\Omega^{k}_{X_{\Gamma}} \left( k \sum_{j} s_{j} \right) \right) 
= (2k-1)(g-1) + kn, 
$$
and denote by $\left\{ \varphi_{1},..., \varphi_{h_{k}} \right\}$ 
its dual basis of 
$H^{0} \left( X_{\Gamma}, 
\Omega^{k}_{X_{\Gamma}} \left( k \sum_{j} s_{j} \right) \right)$ 
for $\Psi_{g,n;k}$, 
and call it {\it normalized}. 
The following theorem will be proved in Section 6. 
\vspace{2ex}

{\bf Theorem 3.4.} 
\begin{it} 
Let $F(k)$ be the holomorphic function given in the formulas of 
Zograf and of McIntyre-Takhtajan.  
Then 
$$
\mu_{g,n;k} 
\left( \varphi_{1} \wedge \cdots \wedge \varphi_{h_{k}} \right) 
= \pm \frac{F(1)^{d_{k}}}{F(k)} 
(\omega_{1} \wedge \cdots \wedge \omega_{g})^{\otimes d_{k}} 
\otimes \bigotimes_{j=1}^{n} d(z - s_{j})^{\otimes (k^{2}-k)/2}. 
$$
\end{it} 

\begin{center}
{\bf 4. Tate curve} 
\end{center}

4.1. \ 
We review the Tate curve (cf. [Si, T]) which is the generalized elliptic curve 
over the ring ${\mathbb Z}[[q]]$ of integral power series of $q$, 
and is the universal elliptic curve over the ring 
${\mathbb Z}((q)) = {\mathbb Z}[[q]][1/q]$ 
of integral Laurent power series of $q$. 
Put 
$$
s_{k} = \sum_{n=1}^{\infty} \frac{n^{k} q^{n}}{1 - q^{n}}, \ \ 
a_{4}(q) = - 5 s_{3}(q), \ \ 
a_{6}(q) = - \frac{5 s_{3}(q) + 7 s_{5}(q)}{12} 
$$
which belong to ${\mathbb Z}[[q]]$. 
Then the Tate elliptic curve ${\cal E}_{q}$ is defined as 
$$
y^{2} + xy = x^{3} + a_{4}(q) x + a_{6}(q), 
$$
and is formally represented as 
${\mathbb G}_{m} / \langle q \rangle$ with the origin $o$ 
given by $1 \in {\mathbb G}_{m}$. 
Therefore, 
$dz/z$ $(z \in {\mathbb G}_{m})$ is a regular $1$-form on ${\cal E}_{q}$, 
and 
\begin{eqnarray*}
X(z) & = & 
\sum_{n \in {\mathbb Z}} \frac{q^{n} z}{(1 - q^{n} z)^{2}} - 2 s_{1}(q) 
\\ & = & 
\frac{z}{(1 - z)^{2}} + \sum_{n=1}^{\infty} 
\left( \frac{q^{n} z}{(1 - q^{n} z)^{2}} + 
\frac{q^{n} z^{-1}}{(1 - q^{n} z^{-1})^{2}} - 
2 \frac{q^{n}}{(1 - q^{n})^{2}} \right), 
\\ 
Y(z) & = & 
\sum_{n \in {\mathbb Z}} \frac{(q^{n} z)^{2}}{(1 - q^{n} z)^{3}} + s_{1}(q) 
\\ & = & 
\frac{z^{2}}{(1 - z)^{3}} + \sum_{n=1}^{\infty} 
\left( \frac{(q^{n} z)^{2}}{(1 - q^{n} z)^{3}} + 
\frac{q^{n} z^{-1}}{(1 - q^{n} z^{-1})^{3}} + 
\frac{q^{n}}{(1 - q^{n})^{2}} \right) 
\end{eqnarray*}
are meromorphic functions on ${\cal E}_{q}$ which have only one pole at $o$ 
of order $2$ and $3$ respectively. 
Hence for each positive integer $k$, 
$H^{0} \left( {\cal E}_{q}, \Omega^{k}_{{\cal E}_{q}} (ko) \right)$ 
is a free ${\mathbb Z}((q))$-module of rank $k$ generated by 
$(dz/z)^{k}$ and $\alpha_{i}(z) (dz/z)^{k}$ $(2 \leq i \leq k)$, 
where $\alpha_{i}(z)$ is a meromorphic function on ${\cal E}_{q}$ 
with only one pole at $o$ of order $i$. 
\vspace{2ex}

4.2. \ 
We consider the Tate curve with marked points. 
For variables $t_{1},..., t_{h}$, 
put 
$$
R = {\mathbb Z} \left[ t_{k}, \frac{1}{t_{k}}, \frac{1}{t_{k}-1}, 
\frac{1}{t_{l} - t_{m}} \ (1 \leq k, l, m \leq h, \ l \neq m) \right]. 
$$
Then each $t_{i}$ gives a point on ${\cal E}_{q} \otimes R((q))$ 
which we denote by the same symbol, 
and 
$$
V = H^{0} \left( {\cal E}_{q} \otimes R((q)), 
\Omega^{k}_{{\cal E}_{q} \otimes R((q))} 
\left( k \sum_{j=1}^{h} t_{j} \right) \right)
$$
is a free $R((q))$-module of rank $hk$. 
Its basis consists of certain products of 
$X(z/t_{j})/t_{j}^{2}$, $Y(z/t_{j})/t_{j}^{3}$ times $(dz/z)^{k}$. 
If we put $t_{1} =1$, 
then the corresponding point becomes the origin $o$. 
\vspace{2ex} 

{\bf Proposition 4.1.} 
\begin{it} 
There exists an $R((q))$-basis of the above $V$ 
which gives the normalized basis. 
\end{it}
\vspace{2ex}

{\it Proof.} 
This follows from that ${\rm Res}_{z=0} \left( z^{k-1} (dz/z^{k}) \right) = 1$ 
and that $X(z/t_{j})/t_{j}^{2}$, $Y(z/t_{j})/t_{j}^{3}$ have Laurent expansion 
at $z = 0$ with coefficients in $R((q))$. 
\ $\square$ 
\vspace{2ex}

\begin{center}
{\bf 5. Arithmetic Schottky uniformization} 
\end{center}

5.1. \ 
Arithmetic Schottky uniformization theory [I2] constructs 
a higher genus version of the Tate curve, 
and its $1$-forms and periods. 
We review this theory for the special case concerned with 
universal deformations of irreducible degenerate curves. 

Denote by $\Delta = \Delta_{g}$ the graph 
with one vertex and $g$ loops. 
Let $x_{\pm 1},..., x_{\pm g}$, $y_{1},..., y_{g}$ be variables, 
and put 
\begin{eqnarray*}
A_{g} & = & {\mathbb Z} \left[ x_{k}, \frac{1}{x_{l}-x_{m}} \ 
\left( k, l, m \in \{ \pm 1,..., \pm g \}, \ l \neq m \right) \right], 
\\
A_{\Delta} & = & A_{g} [[ y_{1},... y_{g} ]], 
\\
B_{\Delta} & = & A_{\Delta} \left[ 1/y_{i} \ (1 \leq i \leq g) \right]. 
\end{eqnarray*}
Then it is shown in [I2, Section 3] that 
there exists a stable curve $C_{\Delta}$ of genus $g$ over $A_{\Delta}$ 
which satisfies the followings: 
\begin{itemize}

\item 
$C_{\Delta}$ is a universal deformation of the universal degenerate curve 
with dual graph $\Delta$ which is obtained from ${\mathbb P}^{1}_{A_{g}}$ 
by identifying $x_{i}$ and $x_{-i}$ $(1 \leq i \leq g)$. 
The ideal of $A_{\Delta}$ generated by $y_{1},..., y_{g}$ 
corresponds to the closed substack 
$\partial {\cal M}_{g} = \overline{\cal M}_{g} - {\cal M}_{g}$ 
of $\overline{\cal M}_{g}$ via the morphism 
${\rm Spec}(A_{\Delta}) \rightarrow \overline{\cal M}_{g}$ 
associated with $C_{\Delta}$. 

\item 
$C_{\Delta}$ is smooth over $B_{\Delta}$, 
and is Mumford uniformized (cf. [M1]) by the subgroup $\Gamma_{\Delta}$ of 
$PGL_{2} \left( B_{\Delta} \right)$ 
with $g$ generators 
$$
\phi_{i} = 
\left( \begin{array}{cc} x_{i} & x_{-i} \\ 1 & 1     \end{array} \right) 
\left( \begin{array}{cc} 1     & 0      \\ 0 & y_{i} \end{array} \right) 
\left( \begin{array}{cc} x_{i} & x_{-i} \\ 1 & 1     \end{array} \right)^{-1} 
\ {\rm mod} \left( B_{\Delta}^{\times} \right) 
\ (1 \leq i \leq g). 
$$
Furthermore, 
$C_{\Delta}$ has the following universality: 
for a complete integrally closed noetherian local ring $R$ 
with quotient field $K$ 
and a Mumford curve $C$ over $K$ such that $\Delta$ is the dual graph of 
its degenerate reduction, 
there is a ring homomorphism $A_{\Delta} \rightarrow R$ 
which gives rise to $C_{\Delta} \otimes_{A_{\Delta}} K \cong C$. 

\item 
By substituting $\alpha_{\pm i} \in {\mathbb C}$ to $x_{\pm i}$ and 
$q_{i} \in {\mathbb C}^{\times}$ to $y_{i}$ $(1 \leq i \leq g)$, 
$C_{\Delta}$ becomes the Riemann surface Schottky uniformized by 
$\Gamma = \langle \gamma_{1},..., \gamma_{g} \rangle$ 
if $\alpha_{\pm i}$ are mutually different and $q_{i}$ are sufficiently small. 

\end{itemize}

Actually, $C_{\Delta}$ is constructed in [I2] as the quotient of 
a certain subspace of ${\mathbb P}^{1}_{B_{\Delta}}$ by the action of 
$\Gamma$ using the theory of formal schemes. 
Furthermore, as is shown in [MD] and [I1, Section 3], 
the normalized holomorphic $1$-forms $\omega_{i}$ $(1 \leq i \leq g)$ 
on Schottky uniformized Riemann surfaces $X_{\Gamma}$ have 
the universal expression 
\begin{eqnarray*}
& & 
\sum_{\phi \in \Gamma_{\Delta} / \langle \phi_{i} \rangle} 
\left( \frac{1}{z - \phi(x_{i})} - \frac{1}{z - \phi(x_{-i})} \right) 
\\
& = & 
\left( \frac{1}{z - x_{i}} - \frac{1}{z - x_{-i}} \right) + 
\sum_{\phi \in \Gamma_{\Delta} / \langle \phi_{i} \rangle - \{ 1 \}} 
\left( \frac{1}{z - \phi(x_{i})} - \frac{1}{z - \phi(x_{-i})} \right) + 
\cdots 
\\
& \in & 
A_{\Delta} \left[ \prod_{k=1}^{g} \frac{1}{(z - x_{k})(z - x_{-k})} \right] 
\end{eqnarray*}
which make a basis of regular $1$-forms on $C_{\Delta}$. 
We denote this basis on $C_{\Delta}$ by the same symbol $\{ \omega_{i} \}$. 

Let $\Delta_{g-1}$ be the graph with one vertex and $(g-1)$ loops, 
and put 
$$
A_{g-1} = {\mathbb Z} \left[ x_{k}, \frac{1}{x_{l}-x_{m}} \ 
\left( k, l, m \in \{ \pm 1,..., \pm (g-1) \}, \ l \neq m \right) \right]. 
$$
Then we have the generalized Tate curve $C_{\Delta_{g-1}}$ of genus $g-1$ 
over $A_{g-1}[[ y_{1},..., y_{g-1}]]$. 
\vspace{2ex} 

{\bf Proposition 5.1.} 
\begin{it} 

{\rm (1)} 
The stable curve $C_{\Delta}|_{y_{g} = 0}$ is obtained from 
$C_{\Delta_{g-1}} \otimes_{A_{g-1}} A_{g}$ by identifying $x_{g} = x_{-g}$. 

{\rm (2)} 
The set $\left\{ \omega_{i}|_{y_{g} = 0} \ | \ 1 \leq i \leq g \right\}$ 
is the normalized basis of 
$H^{0} \left( C_{\Delta}|_{y_{g}=0}, \omega_{C_{\Delta}|_{y_{g}=0}} \right)$. 

\end{it} 
\vspace{2ex}

{\it Proof.} 
The assertion (1) follows from the universality of 
$C_{\Delta}$ and $C_{\Delta_{g-1}}$. 
The assertion (2) follows from that 
$\omega_{i}|_{y_{g} = 0}$ $(i < g)$ are normalized regular $1$-forms 
on $C_{\Delta_{g-1}}$ and that 
$$
\omega_{g}|_{y_{g} = 0} = 
\left( \frac{1}{z - x_{g}} - \frac{1}{z - x_{-g}} + \cdots \right) dz. 
$$
These facts are implied by the variational formula (cf. [F, Chapter III]), 
and also by the universal expression of $\omega_{i}$. 
\ $\square$ 
\vspace{2ex}

5.2. \ 
In what follows, $\phi_{1},..., \phi_{g}$ are {\it normalized} 
by considering $x_{1}, x_{-1}$ as $0, \infty$ respectively, namely, 
$$
\phi_{1} = 
\left( \begin{array}{cc} 0 & 1 \\ 1 & 0     \end{array} \right) 
\left( \begin{array}{cc} 1 & 0 \\ 0 & y_{1} \end{array} \right) 
\left( \begin{array}{cc} 0 & 1 \\ 1 & 0     \end{array} \right)^{-1} 
\ {\rm mod} \left( B_{\Delta}^{\times} \right), 
$$
and by putting $x_{2} = 1$. 
Then as is shown in [I2, 1.1], 
the associated generalized Tate curve $C_{\Delta}$ 
is defined over $A'_{\Delta} = A'_{g} [[y_{1},..., y_{g}]]$, 
where $A'_{g}$ is obtained from $A_{g}$ by deleting $x_{-1}$ and 
putting $x_{1} = 0$, $x_{2} = 1$. 
\vspace{2ex} 

{\bf Proposition 5.2.} 
\begin{it} 

{\rm (1)} 
The infinite products 
$$
\prod_{\{ \gamma \}} \prod_{m=0}^{\infty} \left( 1 - q_{\gamma}^{1+m} \right), 
\ \ 
\left( 1 - q_{\gamma_{1}} \right)^{2} \cdots 
\left( 1 - q_{\gamma_{1}}^{k-1} \right)^{2} 
\left( 1 - q_{\gamma_{2}}^{k-1} \right) 
\prod_{\{ \gamma \}} \prod_{m=0}^{\infty} \left( 1 - q_{\gamma}^{k+m} \right) 
$$
given in 2.2 have universal expression as elements of $A'_{\Delta}$. 

{\rm (2)} 
Under $y_{2} = \cdots = y_{g} = 0$, 
these elements of $A'_{\Delta}$ becomes 
$\prod_{m=0}^{\infty} \left( 1 - y_{1}^{1+m} \right)^{2}$ 
which is primitive, 
i.e., not congruent to $0$ modulo any rational prime. 
\end{it}
\vspace{2ex}

{\it Proof.} 
Let $(\Gamma; \gamma_{1},..., \gamma_{g})$ be a normalized Schottky group, 
and for $i = 1,..., g$, put $\gamma_{-i} = \gamma_{i}^{-1}$. 
Then by Proposition 1.3 of [I2] and its proof, 
if $\gamma \in \Gamma$ has the reduced expression 
$\gamma_{\sigma(1)} \cdots \gamma_{\sigma(l)}$ 
$(\sigma(i) \in \{ \pm 1,..., \pm g \})$ 
such that $\sigma(1) \neq - \sigma(l)$, 
then its multiplier $q_{\gamma}$ has universal expression as an element 
of $A'_{\Delta}$ divisible by $y_{\sigma(1)} \cdots y_{\sigma(l)}$. 
Therefore, the assertion (1) holds. 
Since 
$$
\left( 1 - y_{1} \right)^{2} \cdots \left( 1 - y_{1}^{k-1} \right)^{2} 
\prod_{m=0}^{\infty} \left( 1 - y_{1}^{k+m} \right)^{2}  
= \prod_{m=0}^{\infty} \left( 1 - y_{1}^{1+m} \right)^{2}, 
$$
the assertion (2) holds. 
\ $\square$ 
\vspace{2ex}

\begin{center}
{\bf 6. Product formula of Mumford forms} 
\end{center}

6.1. \ 
First, we will prove the following: 
\vspace{2ex}

{\bf Theorem 6.1.} 
\begin{it} 
If $n = 0$, then Theorem 3.4 holds. 
\end{it}
\vspace{2ex} 

{\it Proof.} 
As in Section 5, 
let $C_{\Delta}$ be the generalized Tate curve over $A'_{\Delta}$ of  genus $g$ 
which is uniformized by $\Gamma_{\Delta}$, 
where the generators $\phi_{1},..., \phi_{g}$ of $\Gamma_{\Delta}$ 
is normalized as $x_{1} = 0$, $x_{-1} = \infty$, $x_{2} = 1$. 
As in stated in 5.1, 
we can consider $C_{\Delta}$ as a family of Schottky uniformized 
Riemann surfaces by taking $x_{\pm i}$, $y_{i}$ 
as complex parameters $\alpha_{\pm i}$, $q_{i}$ respectively 
such that $q_{i}$ are sufficiently small. 
By Proposition 5.1, 
$C_{\Delta}|_{y_{2} = \cdots = y_{g} = 0}$ becomes the stable curve 
$C'_{\Delta}$ obtained from the Tate curve 
${\cal E}_{y_{1}} = {\mathbb G}_{m}/\langle y_{1} \rangle$ 
over $A'_{g} [[y_{1}]]$ by identifying $x_{i} = x_{-i}$ $(2 \leq i \leq g)$. 
Hence by the properties of the relative dualizing sheaf [K, Section 1], 
$H^{0} \left( C'_{\Delta}, \omega^{k}_{C'_{\Delta}} \right)$ 
becomes the subspace of 
$$
W = H^{0} \left( {\cal E}_{y_{1}}, \omega^{k}_{{\cal E}_{y_{1}}} 
\left( k \sum_{i=2}^{g} (x_{i} + x_{-i}) \right) \right)
$$ 
which consists of $\eta \in W$ satisfying 
${\rm Res}^{k}_{x_{i}} (\eta) = (-1)^{k} {\rm Res}^{k}_{x_{-i}} (\eta)$. 
By Proposition 4.1, 
there exists the normalized basis of $W$ as an $A'_{g}((y_{1}))$-basis, 
and hence there exists a basis $\{ \phi_{l} \}$ of 
$H^{0} \left( C_{\Delta}, \omega^{k}_{C_{\Delta}} \right) \otimes 
A'_{\Delta}[1/y_{1}]$ 
which converges to the normalized basis 
if $y_{2},..., y_{g} \rightarrow 0$. 
For each $i = 2,..., g$, 
let 
$$
z_{i} = \frac{(x_{i} - x_{-i})(z - x_{i})}{z - x_{-i}}, \   
z_{-i} = \frac{(x_{-i} - x_{i})(z - x_{-i})}{z - x_{i}} 
$$
be the local coordinates at $x_{i}$, $x_{-i}$ respectively 
such that 
$$
\lim_{z \rightarrow x_{i}} (z - x_{i})/z_{i} = 
\lim_{z \rightarrow x_{-i}} (z - x_{-i})/z_{i} = 1. 
$$
Then the transformation matrices of 
$$
\left( z^{j} \right)_{0 \leq j \leq l} \mapsto 
\left( (z - x_{\pm i})^{j} \right)_{0 \leq j \leq l}, \ 
\left( (z - x_{\pm i})^{j} \right)_{0 \leq j \leq l} \mapsto 
\left( z_{\pm i}^{j} \ {\rm mod} \left( z_{\pm i}^{l+1} \right) 
\right)_{0 \leq j \leq l} 
$$
have determinant $1$, 
and hence we may replace the exterior product of 
the normalized basis with that of the basis of 
$H^{0} \left( C_{\Delta}, \omega^{k}_{C_{\Delta}} \right)$ 
dual to $\left( z_{\pm i}^{j} \right)_{1 \leq i \leq g}$. 
Since $C_{\Delta}$ is obtained from ${\cal E}_{y_{1}}$ by the equation 
$z_{i} z_{-i} = -(x_{i} - x_{-i})^{2} y_{i}$ $(2 \leq i \leq g)$, 
$$
\frac{(dz_{i})^{k}}{z_{i}^{k+1}} \wedge \cdots \wedge 
\frac{(dz_{i})^{k}}{z_{i}^{2k-1}} 
= \pm \left( (x_{i} - x_{-i})^{2} y_{i} \right)^{-(k^{2}-k)/2} 
\frac{(dz_{-i})^{k}}{z_{-i}} \wedge \cdots \wedge 
\frac{(dz_{-i})^{k}}{z_{-i}^{k-1}}. 
$$
Further, 
$\frac{1}{2 \pi \sqrt{-1}} \oint_{C_{\pm i}} \rightarrow 
{\rm Res}_{x_{\pm i}}$ 
under $y_{i} \rightarrow 0$. 
Hence if $k > 1$ and $y_{2},..., y_{g} \rightarrow 0$, 
then 
$$ 
\prod_{i=2}^{g} \left( (x_{i} - x_{-i})^{2} y_{i} \right)^{(k^{2}-k)/2} 
\det \left( \Psi_{g;k} (\phi_{l}, \xi_{m}) \right)_{l,m} \rightarrow \pm 1, 
$$
where $\{ \xi_{m} \}$ is the natural basis 
$$
\left\{ \xi_{1,k-1}, \xi_{2,1},..., \xi_{2,2k-2}, 
\xi_{i,0},..., \xi_{i,2k-2} \ (3 \leq i \leq g) \right\}
$$
given in 2.2. 
Therefore, 
for the normalized basis $\{ \varphi_{l} \}$ of holomorphic $k$-forms 
on $C_{\Delta}$, 
$\prod_{i=2}^{g} \left( (x_{i} - x_{-i})^{2} y_{i} \right)^{(k^{2}-k)/2} 
\bigwedge_{l} \varphi_{l}$ 
becomes the exterior product of the normalized basis of $W$ up to a sign 
under $y_{2},..., y_{g} \rightarrow 0$. 
Hence by Theorem 3.1 (in the case when $n = 0$), 
Propositions 3.2, 5.1 and 5.2, 
if $y_{2},..., y_{g} \rightarrow 0$, 
then 
$$
c(g;k) \cdot \mu_{1,2g-2;k} = 
\pm \prod_{m=1}^{\infty} \left( 1 - y_{1}^{m} \right)^{24} 
\times \prod_{i=2}^{g} (x_{i} - x_{-i})^{2} 
$$
under the evaluation on the normalized basis of the above $W$. 
Therefore, by Proposition 4.1, 
we have $c(g;k) = \pm 1$. 
\ $\square$ 
\vspace{2ex}

{\bf Corollary 6.2.} 
\begin{it} 
Let $c_{g}$ and $c_{g;k}$ be the constants given in the formulas 
of Zograf and of McIntyre-Takhtajan respectively. 
Then 
$$
c_{g;k} = (c_{g})^{d_{k}} 
\exp \left[ (1-g)(k^{2}-k) \left( 24 \zeta'_{\mathbb Q}(-1)-1 \right) \right]. 
$$
\end{it}
\vspace{-2ex}

{\it Proof.} 
This follows from Proposition 2.1 and Theorem 6.1. 
\ $\square$ 
\vspace{2ex} 

6.2. \ 
We give a sketch of another proof of Theorem 6.1 in the case where $k = 2$. 
Let $\zeta_{i,j}$ $(1 \leq i \leq g, \ 0 \leq j \leq 2)$ be the map from 
$\left\{ \phi_{1},..., \phi_{g} \right\}$ into $A'_{g}[z]_{2}$ defined as 
$$
\zeta_{i,j} \left( \phi_{l} \right) = \delta_{il} (z - x_{i})^{j} 
\ (1 \leq l \leq g). 
$$
By a theorem of Max Noether, 
the space of regular $2$-forms on a genus $2$ or non-hyperelliptic curve 
is generated by the products of regular $1$-forms. 
Since the paring given by $\Psi$ between $\omega_{k} \omega_{l}$ 
and $\zeta_{i,j}$ has universal expression as an element of $A'_{\Delta}$, 
to prove the assertion, 
it is enough to show that for a field $K$ of any characteristic, 
the matrix consisting of the pairings between 
$$
\left\{ \omega_{l}^{2} \ (1 \leq l \leq g), 
\ \omega_{1} \omega_{l} \ (2 \leq l \leq g),  
\ \omega_{2} \omega_{l} \ (3 \leq l \leq g) \right\} 
$$
and 
$$
\left\{ \zeta_{1,1}, \ \zeta_{2,1}, \ \zeta_{2,2}, \ 
\zeta_{i,0}, \ \zeta_{i,1}, \ \zeta_{i,2} \ (3 \leq i \leq g) \right\} 
$$
is invertible over the quotient field $K_{\Delta}$ of 
$A'_{\Delta} \widehat{\otimes}_{\mathbb Z} K$. 

Let $I_{\Delta}$ be the ideal of $A'_{\Delta}$ 
generated by $y_{1},..., y_{g}$. 
Then we have 
$$
\left\{ \begin{array}{ll} 
\Psi \left( \omega_{1} \omega_{l}, \zeta_{i,0} \right) \equiv 
{\displaystyle \frac{\delta_{il}(x_{1} - x_{-1})}
{(x_{i} - x_{1})(x_{i} - x_{-1})}} \ {\rm mod} \left( I_{\Delta} \right) 
& (2 \leq l \leq g), 
\\ 
\Psi \left( \omega_{2} \omega_{l}, \zeta_{i,0} \right) \equiv 
{\displaystyle \frac{\delta_{il}(x_{2} - x_{-2})}
{(x_{i} - x_{2})(x_{i} - x_{-2})}} \ {\rm mod} \left( I_{\Delta} \right) 
& (3 \leq l \leq g) 
\end{array} \right. 
$$
for $3 \leq i \leq g$, 
$$
\left\{ \begin{array}{ll} 
\Psi \left( \omega_{l}^{2}, \zeta_{i,1} \right) \equiv 
\delta_{il} \ {\rm mod} \left( I_{\Delta} \right) 
& (1 \leq l \leq g), 
\\
\Psi \left( \omega_{1} \omega_{l}, \zeta_{i,1} \right) \equiv 
0 \ {\rm mod} \left( I_{\Delta} \right) 
& (2 \leq l \leq g), 
\\
\Psi \left( \omega_{2} \omega_{l}, \zeta_{i,1} \right) \equiv 
0 \ {\rm mod} \left( I_{\Delta} \right) 
& (3 \leq l \leq g) 
\end{array} \right. 
$$
for $1 \leq i \leq g$, 
and 
$$
\left\{ \begin{array}{ll} 
\Psi \left( \omega_{1} \omega_{l}, \zeta_{i,2} \right) \equiv 
{\displaystyle 
\frac{\delta_{il}(x_{1} - x_{-1})(x_{i} - x_{-i})^{2}}
{(x_{-i} - x_{1})(x_{-i} - x_{-1})} y_{i}} 
\ {\rm mod} \left( I_{\Delta}^{2} \right) 
& (2 \leq l \leq g), 
\\
\Psi \left( \omega_{2} \omega_{l}, \zeta_{2,2} \right) \equiv 
{\displaystyle 
\frac{(x_{l} - x_{-l})(x_{2} - x_{-2})^{2}}
{(x_{l} - x_{-2})(x_{-l} - x_{-2})} y_{2}} 
\ {\rm mod} \left( I_{\Delta}^{2} \right), 
& 
\\
\Psi \left( \omega_{2} \omega_{l}, \zeta_{i,2} \right) \equiv 
{\displaystyle 
\frac{\delta_{il}(x_{2} - x_{-2})(x_{i} - x_{-i})^{2}}
{(x_{-i} - x_{2})(x_{-i} - x_{-2})} y_{i}} 
\ {\rm mod} \left( I_{\Delta}^{2} \right) 
& (i \neq 2, \ 3 \leq l \leq g) 
\end{array} \right. 
$$
for $2 \leq i \leq g$. 
Therefore, 
the lowest term of this determinant is the product 
$$
\prod_{k=3}^{g} \sigma_{k} \cdot \prod_{k=2}^{g} \tau_{k}, 
$$
where 
$$
\sigma_{k} = 
\frac{x_{2} - x_{-2}}{(x_{k} - x_{2})(x_{k} - x_{-2})}, \ 
\tau_{2} = 
\frac{(x_{1} - x_{-1})(x_{2} - x_{-2})^{2}}
{(x_{-2} - x_{1})(x_{-2} - x_{-1})} y_{2}, 
$$
and for $3 \leq k \leq g$, 
\begin{eqnarray*}
\tau_{k} 
& \equiv & 
\Psi \left( \omega_{1} \omega_{k}, \zeta_{k,2} \right) - 
\Psi \left( \omega_{2} \omega_{k}, \zeta_{k,2} \right) 
\frac{\Psi \left( \omega_{1} \omega_{k}, \zeta_{k,0} \right)}
{\Psi \left( \omega_{2} \omega_{k}, \zeta_{k,0} \right)} 
\ {\rm mod} \left( I_{\Delta}^{2} \right) 
\\
& = & 
\left( 
\frac{1}{(x_{-k} - x_{1})(x_{-k} - x_{-1})} - 
\frac{(x_{k} - x_{2})(x_{k} - x_{-2})}
{(x_{-k} - x_{2})(x_{-k} - x_{-2})(x_{k} - x_{1})(x_{k} - x_{-1})} 
\right) 
\\
& & \times \ (x_{1} - x_{-1})(x_{k} - x_{-k})^{2} y_{k}. 
\end{eqnarray*}
Hence this product is invertible in $K_{\Delta}$ and the assertion holds. 
\vspace{2ex}

6.3. \ 
Second, we deduce Theorem 3.4 from Theorem 6.1. 
\vspace{2ex}

{\it Proof of Theorem 3.4.} 
Let $\left( \pi : C \rightarrow S; \sigma_{1},..., \sigma_{n} \right)$ 
be an $n$-pointed stable curve of genus $g \geq 2$, 
Then for non-negative integers $l_{1},..., l_{n} \leq k$, 
the $l_{i}$th residue map gives the exact sequence 
$$
0 \rightarrow 
\pi_{*} \left( \omega_{C}^{k} \left( \sum_{j=1}^{n} l_{j} \sigma_{j} \right) 
\right) 
\rightarrow 
\pi_{*} \left( \omega_{C}^{k} \left( \sum_{j \neq i} l_{j} \sigma_{j} + 
(l_{i} + 1) \sigma_{i} \right) \right) 
\rightarrow 
{\cal O}_{S} \cdot (d z_{\sigma_{i}})^{k-l_{i}-1} 
\rightarrow 0, 
$$
where $d z_{\sigma_{i}}$ denotes the local coordinate around $\sigma_{i}$. 
Therefore, by induction, 
$$
\det \pi_{*} \left( \omega_{C}^{k} \left( k \sum_{j=1}^{n} \sigma_{j} 
\right) \right) 
\cong 
\det \pi_{*} \left( \omega_{C}^{k} \right) \cdot 
\prod_{j=1}^{n} (d z_{\sigma_{i}})^{(k^{2} - k)/2}. 
$$
If $C$ consists of Schottky uniformized Riemann surfaces 
with marked points $s_{1},..., s_{n}$, 
then the residue map sends the normalized basis to 
$\left\{ d(z - s_{i})^{k - l_{i} - 1} \right\}$, 
and hence the assertion follows from Theorem 6.1. 
\ $\square$ 
\vspace{2ex}

\begin{center}
{\bf 7. Rationality of Ruelle zeta values} 
\end{center}

7.1. \ 
Assume that $g > 1$, 
let $\left( \Gamma; \gamma_{1},..., \gamma_{g} \right)$ be 
a marked normalized Schottky group, 
and $F_{\left( \Gamma; \gamma_{1},..., \gamma_{g} \right)}(k)$ 
be the value of $F(k)$ at the corresponding point on ${\mathfrak S}_{g}$. 
Let $R$ be a subring of ${\mathbb C}$ containing $1$ such that 
$X_{\Gamma}$ has a model ${\mathcal X}$ as a stable curve 
over $R$, 
and assume that there exist $R$-bases $\{ u_{1},..., u_{g} \}$ of 
$H^{0} \left( {\cal X}, \omega_{{\cal X}/R} \right)$ 
and $\{ v_{1},..., v_{(2k-1)(g-1)} \}$ of 
$H^{0} \left( {\cal X}, \omega^{k}_{{\cal X}/R} \right)$ for $k > 1$. 
Then their periods are defined as 
$$
\Omega_{1} = \det \left( \frac{1}{2 \pi \sqrt{-1}} 
\oint_{C_{i}} u_{j} \right)_{1 \leq i, j \leq g}, 
$$ 
and  
$$
\Omega_{k} = \det \left( \frac{1}{2 \pi \sqrt{-1}} \sum_{i=1}^{g} 
\oint_{C_{i}} v_{l} \cdot \xi_{m}(\gamma_{i}) 
\right)_{1 \leq l,m \leq (2k-1)(g-1)}, 
$$
where 
$\left\{ \xi_{1},..., \xi_{(2k-1)(g-1)} \right\}$ is the basis 
$\left\{ \xi_{1,k-1}, \ \xi_{2,1},..., \xi_{2,2k-2}, \ 
\xi_{i,0},..., \xi_{i,2k-2} \ (3 \leq i \leq g) \right\}$ 
of $H^{1} \left( \Gamma, {\mathbb C}[z]_{2k-2} \right)$ 
given in 2.2. 
We consider the rationality of  
$$
\frac{F_{\left( \Gamma; \gamma_{1},..., \gamma_{g} \right)}(k)}{\Omega_{k}}, 
$$
and its relation with the discriminant ideal $D_{\cal X}$ 
of ${\cal X}$. 
Here $D_{\cal X}$ is defined as an ideal of $R$ which corresponds 
to the closed substack 
$\partial {\cal M}_{g} = \overline{\cal M}_{g} - {\cal M}_{g}$ 
of $\overline{\cal M}_{g}$ via the morphism 
${\rm Spec} (R) \rightarrow \overline{\cal M}_{g}$ 
associated with ${\cal X}$. 
\vspace{2ex}

{\bf Theorem 7.1.} 
\begin{it} 

{\rm (1)} 
Under the above notations,  
the ratio 
$$
\frac{\left( F_{\left( \Gamma; \gamma_{1},..., \gamma_{g} \right)}(1)
/\Omega_{1} \right)^{d_{k}}}
{F_{\left( \Gamma; \gamma_{1},..., \gamma_{g} \right)}(k)
/\Omega_{k}} 
$$
belongs to the power $D_{\cal X}^{(k^{2}-k)/2}$ of $D_{\cal X}$. 

{\rm (2)} 
If $R$ is a discrete valuation ring and 
$d_{\cal X}$ is a generator of $D_{\cal X}$, 
then 
$$
\frac{\left( F_{\left( \Gamma; \gamma_{1},..., \gamma_{g} \right)}(1)
/\Omega_{1} \right)^{d_{k}}}
{F_{\left( \Gamma; \gamma_{1},..., \gamma_{g} \right)}(k)
/\Omega_{k}} 
\in d_{\cal X}^{(k^{2}-k)/2} \cdot R^{\times}, 
$$
where $R^{\times}$ denotes the unit group of $R$. 
 
{\rm (3)} 
If ${\cal X}$ is smooth over $R$, 
then 
$$
\frac{\left( F_{\left( \Gamma; \gamma_{1},..., \gamma_{g} \right)}(1)
/\Omega_{1} \right)^{d_{k}}}
{F_{\left( \Gamma; \gamma_{1},..., \gamma_{g} \right)}(k)
/\Omega_{k}} 
\in R^{\times}.
$$
\end{it}

{\it Proof.} 
By Theorem 6.1, 
this ratio is (up to a sign) the evaluation of the Mumford form  $\mu_{g;k}$ on ${\mathcal X}$ under the trivializations of $\lambda_{1}$ and $\lambda_{k}$ 
by $u_{1} \wedge \cdots \wedge u_{g}$ and 
$v_{1} \wedge \cdots \wedge v_{(2k-1)(g-1)}$ respectively. 
Hence the assertions follows from by Theorem 3.1. 
\ $\square$ 
\vspace{2ex}

7.2. \ 
Let $\Gamma$ be a Schottky group, 
and denote by ${\mathbb H}^{3}/\Gamma$ the hyperbolic $3$-fold 
uniformized by $\Gamma$. 
Then for each $\gamma \in \Gamma$, 
$-\log |q_{\gamma}|$ is the length of the closed geodesic 
on ${\mathbb H}^{3}/\Gamma$ corresponding to $\gamma$, 
and hence the Ruelle zeta function of ${\mathbb H}^{3}/\Gamma$ becomes 
$$
R_{\Gamma}(s) = \prod_{\{ \gamma \}} \left( 1 - |q_{\gamma}|^{s} \right)^{-1}, 
$$
where $\{ \gamma \}$ runs over primitive conjugacy classes of $\Gamma$. 
It is known (cf. [MT, 5.2]) that $R_{\Gamma}(s)$ is absolutely convergent 
if ${\rm Re}(s) \geq 2$. 
We assume that a marked normalized Schottky group 
$( \Gamma; \gamma_{1},..., \gamma_{g})$ is contained in 
$PSL_{2}({\mathbb R})$, 
and apply Theorem 7.1 to showing the rationality of 
the {\it modified} Ruelle zeta values 
$$
\widetilde{R}_{\Gamma}(k) = R_{\Gamma}(k) 
\frac{\left( 1 - q_{\gamma_{1}}^{k} \right)^{2} 
\left( 1 - q_{\gamma_{2}}^{k} \right)}
{\left( 1 - q_{\gamma_{2}}^{k-1} \right)} 
$$
for integers $k > 1$. 
\vspace{2ex}

{\bf Theorem 7.2.} 
\begin{it} 
Assume that $\Gamma \subset PSL_{2}({\mathbb R})$ and 
$k$ is an integer $> 1$. 
Put 
$$
c(\Gamma) = \frac{F_{(\Gamma; \gamma_{1},..., \gamma_{g})}(1)}{\Omega_{1}}. 
$$ 

{\rm (1)} 
If $K$ denotes the quotient field of $R$, 
then 
$$
\widetilde{R}_{\Gamma}(k) 
\in 
\frac{\Omega_{k+1}}{\Omega_{k}} \cdot c(\Gamma)^{12k} \cdot K^{\times}. 
$$

{\rm (2)} 
If $R$ is a discrete valuation ring, 
then 
$$
\widetilde{R}_{\Gamma}(k) 
\in 
\frac{\Omega_{k+1}}{\Omega_{k}} \cdot c(\Gamma)^{12k} \cdot d_{\chi}^{-k} 
\cdot R^{\times}. 
$$

{\rm (3)} 
If ${\cal X}$ is smooth over $R$, 
then 
$$
\widetilde{R}_{\Gamma}(k) 
\in 
\frac{\Omega_{k+1}}{\Omega_{k}} \cdot c(\Gamma)^{12k} \cdot R^{\times}. 
$$
\end{it}

{\it Proof.} 
By the assumption, 
$q_{\gamma} = |q_{\gamma}|$ for any $\gamma \in \Gamma - \{ 1 \}$, 
and hence 
$$
\widetilde{R}_{\Gamma}(k) = 
\frac{F_{(\Gamma; \gamma_{1},..., \gamma_{g})}(k+1)}
{F_{(\Gamma; \gamma_{1},..., \gamma_{g})}(k)}. 
$$
Therefore, 
the assertion follows from Theorem 7.1. 
\ $\square$ 
\vspace{2ex}

{\bf Corollary 7.3.} 
\begin{it} 
Assume that a Schottky group $\Gamma$ is contained in $PSL_{2}({\mathbb R})$, 
and a subring $R$ of ${\mathbb C}$ is a Dedekind ring such that 
$X_{\Gamma}$ has a model ${\cal X}$ as a stable curve over $R$. 
Denote by $K$ the quotient field of $R$. 
For $j = 1, k, k+1$, 
let $\left\{ v_{i}^{(j)} \right\}_{i}$ be a $K$-basis of 
$H^{0} \left( {\cal X} \otimes K, \omega^{j}_{{\cal X} \otimes K/K} \right)$, 
and denote by $\Omega_{j}$ the period of this basis. 

{\rm (1)} 
Put 
$$
c(\Gamma) = \frac{F_{(\Gamma; \gamma_{1},..., \gamma_{g})}(1)}{\Omega_{1}}. 
$$ 
Then there exists an element $r_{\cal X}$ of $K^{\times}$ such that 
$$
\widetilde{R}_{\Gamma}(k) 
\in 
\frac{\Omega_{k+1}}{\Omega_{k}} \cdot c(\Gamma)^{12k} \cdot r_{\cal X}. 
$$

{\rm (2)} 
For each prime ideal $P$ of $R$, 
${\rm ord}_{P} \left( r_{\cal X} \right)$ is given by 
$$
k \cdot {\rm ord}_{P} \left( d_{\cal X} \right) 
+ 12 k \cdot {\rm ord}_{P} \left( \bigwedge v_{i}^{(1)} \right) 
+ {\rm ord}_{P} \left( \bigwedge v_{i}^{(k)} \right) 
- {\rm ord}_{P} \left( \bigwedge v_{i}^{(k+1)} \right), 
$$
where ${\rm ord}_{P} \left( \bigwedge v_{i}^{(j)} \right)$ denotes 
the order at $P$ of the fractional ideal of $R_{P}$ generated by 
the exterior product of $\left\{ v_{i}^{(j)} \right\}_{i}$ 
under the identification 
$$
R_{P} = \det \left( H^{0} \left( {\cal X} \otimes R_{P}, 
\omega^{j}_{{\cal X} \otimes R_{P}/R_{P}} \right) \right). 
$$
\end{it}

\begin{center}
{\bf References} 
\end{center}

\begin{itemize}

\item[{[B]}] 
A. Beilinson, 
Higher regulators and values of $L$-functions, 
J. Sov. Math. {\bf 30} (1985) 2036--2070. 

\item[{[BK]}] 
S. Bloch and K. Kato, 
$L$-functions and Tamagawa numbers of motives, 
in: The Grothendieck Festschrift, vol. I, 
Progr. Math. Vol. 86, Birkh\"{a}user, Boston, 1990, pp. 333--400. 

\item[{[D1]}] 
P. Deligne, 
Valeurs de fonctions $L$ et periodes d'integrales, 
in: Automorphic forms, representations, and $L$-functions, 
Proc. Symp. in Pure Math., Vol. 33, Part 2, 
Amer. Math. Soc., 1979, pp. 313--346. 

\item[{[D2]}] 
P. Deligne, 
Le d\'{e}terminant de la cohomologie, 
in: Current trends in arithmetical algebraic geometry, 
Contemp. Math., Vol. 67, Amer. Math. Soc., 1987, pp. 93--177. 

\item[{[F]}] 
J. D. Fay, 
Theta-functions on Riemann surfaces, 
Lecture Notes in Math., Vol. 352, Springer, 1973. 

\item[{[Fre]}] 
G. Freixas i Montplet, 
An arithmetic Hilbert-Samuel theorem for pointed stable curves, 
J. Eur. Math. Soc. {\bf 14} (2012) 321--351. 

\item[{[Fri]}] 
D. Fried, 
Analytic torsion and closed geodesics on hyperbolic manifolds, 
Invent. Math. {\bf 84} (1986) 523--540. 

\item[{[GP]}] 
Y. Gon and J. Park, 
The zeta functions of Ruelle and Selberg for hyperbolic manifolds with cusps, 
Math. Ann. {\bf 346} (2010) 719--767. 

\item[{[GS]}] 
H. Gillet and C. Soul\'{e}, 
Arithmetic intersection theory, 
Pub. Math. IHES {\bf 72} (1990) 94--174. 

\item[{[I1]}] 
T. Ichikawa, 
$P$-adic theta functions and solutions of the KP hierarchy, 
Comm. Math. Phys. {\bf 176} (1996) 383--399. 

\item[{[I2]}] 
T. Ichikawa, 
Generalized Tate curve and integral Teichm\"{u}ller modular forms, 
Amer. J. Math. {\bf 122} (2000) 1139--1174. 

\item[{[I3]}] 
T. Ichikawa, 
Universal periods of hyperelliptic curves and their applications, 
J. Pure Appl. Algebra {\bf 163} (2001) 277--288. 

\item[{[I4]}] 
T. Ichikawa, 
Klein's amazing formula and arithmetic of Teichm\"{u}ller modular forms, 
preprint. 

\item[{[K]}] 
F. F. Knudsen, 
The projectivity of the moduli space of stable curves II, III, 
Math. Scand. {\bf 52} (1983) 161--199, 200--212. 

\item[{[MD]}]
Y. Manin and V. Drinfeld, 
Periods of $p$-adic Schottky groups, 
J. Reine Angew. Math. {\bf 262/263} (1973) 239--247. 

\item[{[MT]}] 
A. McIntyre and L. A. Takhtajan, 
Holomorphic factorization of determinants of Laplacians on Riemann surfaces 
and a higher genus generalization of Kronecker's first limit formula, 
GAFA Geom. funct. anal. {\bf 16} (2006) 1291--1323. 

\item[{[M1]}] 
D. Mumford, 
An analytic construction of degenerating curves over complete local rings, 
Compositio Math. {\bf 24} (1972) 129--174. 

\item[{[M2]}] 
D. Mumford, 
Stability of projective varieties, 
L'Ens. Math. {\bf 23} (1977) 39--110. 

\item[{[P]}] 
J. Park, 
Analytic torsion and Ruelle zeta functions for hyperbolic manifolds 
with cusps, 
J. Funct. Anal. {\bf 257} (2009) 1713--1758. 

\item[{[S]}] 
T. Saito, 
Conductor, discriminant, and the Noether formula of arithmetic surfaces, 
Duke Math. J. {\bf 57} (1988) 151--173. 

\item[{[Si]}] 
J. H. Silverman, 
Advanced topics in the arithmetic of elliptic curves, 
Graduate Texts in Math. Springer, 1994. 

\item[{[Su1]}] 
K. Sugiyama, 
An analog of the Iwasawa conjecture for a compact hyperbolic threefold, 
J. Reine Angew. Math. {\bf 613} (2007) 35--50. 

\item[{[Su2]}] 
K. Sugiyama, 
A special value of Ruelle $L$-function and the theorem of 
Cheeger and M\"{u}ller, 
arXiv:0803.2079v1 

\item[{[Su3]}] 
K. Sugiyama, 
The Taylor expansion of Ruelle $L$-function at the origin 
and the Borel regulator, 
arXiv:0804.2715v1 

\item[{[Su4]}] 
K. Sugiyama, 
On geometric analogues of the Birch and Swinnerton-Dyer conjecture 
for low dimensional hyperbolic manifolds, 
in: Spectral analysis in geometry and number theory, 
Contemp. Math., Vol. 484, Amer. Math. Soc., 2009. pp. 267--286. 

\item[{[Su5]}] 
K. Sugiyama, 
On geometric analogues of Iwasawa main conjecture for a hyperbolic threefold, 
in: Noncommutativity and singularities, 
Adv. Stud. Pure Math., Vol. 55, Math. Soc. Japan, Tokyo, 2009, pp. 117--135. 

\item[{[T]}] 
J. Tate, 
A review of non-archimedean elliptic functions, 
in: Elliptic Curves, Modular forms, \& Fermat's Last Theorem, 
International Press, Boston, 1995, pp. 162--184. 

\item[{[W]}] 
L. Weng, 
$\Omega$-admissible theory II. 
Deligne pairings over moduli spaces of punctured Riemann surfaces, 
Math. Ann. {\bf 320} (2001) 239--283. 

\item[{[Z1]}] 
P. G. Zograf, 
Liouville action on moduli spaces and uniformization of
degenerate Riemann surfaces, 
Algebra i Analiz {\bf 1} (1989) 136--160 (Russian),
English translation in Leningrad Math. J. {\bf 1} (1990) 941--965. 

\item[{[Z2]}] 
P. G. Zograf, 
Determinants of Laplacians, Liouville action, 
and an analogue of the Dedekind $\eta$-function on Teichm\"{u}ller space, 
Unpublished manuscript (1997). 

\item[{[ZT]}]
P. G. Zograf and L. A. Takhtajan, 
On the uniformization of Riemann surfaces and 
on the Weil-Petersson metric on the Teichm\"{u}ller and 
Schottky spaces, 
Math. Sb. (N. S.) {\bf 132 (174)} (1987) 304--321 (Russian), 
English translation in Math. USSR-Sb. {\bf 60} (1988) 297--313. 

\end{itemize}

{\sc Department of Mathematics, Graduate School of Science and Engineering, 
Saga University, Saga 840-8502, Japan} 

{\it E-mail}: ichikawa@ms.saga-u.ac.jp, ichikawn@cc.saga-u.ac.jp

\end{document}